\documentclass{amsart}
\usepackage{mathrsfs}
\usepackage{amsmath,amscd}
\usepackage{amsmath,amsthm,amssymb,float,multicol}
\usepackage{bibarts}
\usepackage{cite}
\usepackage[all,cmtip]{xy}
\usepackage{hyperref}

\usepackage[pdftex]{graphicx}
\usepackage[usenames]{color}
\usepackage[large]{subfigure}
\usepackage{xypic}

\DeclareMathOperator{\N}{\mathbb{N}}

\def\C{\mathscr{C}}

\def\I{\mathcal{I}}

\def\N{\mathcal{N}}

\def\O{\mathcal{O}}
\def\P{\mathbb{P}}

\def\X{\mathcal{X}}

\def\<{\langle}
\def\>{\rangle}
\def\({\left(}
\def\){\right)}

\newtheorem{theorem}{Theorem}[section]
\newtheorem{lemma}[theorem]{Lemma}

\newtheorem{prop}[theorem]{Proposition}
\newtheorem{cor}[theorem]{Corollary}
\newtheorem{question}[theorem]{Question}
\theoremstyle{definition}
\newtheorem{definition}[theorem]{Definition}
\newtheorem{notation}[theorem]{Notation}
\newtheorem{construction}[theorem]{Construction}

\newtheorem*{ack}{Acknowledgment}

\theoremstyle{remark}

\numberwithin{equation}{section}

\begin{document}
\bibliographystyle{alpha}
\title{Fano Hypersurfaces in Positive Characteristic}
\author{Yi Zhu}
\address{Department of Mathematics, Stony Brook University, Stony Brook, NY, 11794}
\curraddr{}
\email{yzhu@math.sunysb.edu}


\date{}

\dedicatory{}

\begin{abstract}
We prove that a general Fano hypersurface in a projective space over an algebraically closed field of arbitrary characteristic is separably rationally connected.
\end{abstract}

\maketitle

\section{Introduction}

In this paper, we work with varieties over an algebraically closed field $k$ of arbitrary characteristic. 

\begin{definition}[\cite{Kollar} IV.3]Let $X$ be a variety defined over $k$.

A variety $X$ is \emph{rationally connected} if there is a family of irreducible proper rational curves $g: U\rightarrow Y$ and an evaluation morphism $u:U\rightarrow X$ such that the morphism $u^{(2)}:U\times_Y U\rightarrow X\times X$ is dominant.

A variety $X$ is \emph{separably rationally connected} if there exists a proper rational curve $f:\P^1\rightarrow X $ such that the image lies in the smooth locus of $X$ and the pullback of the tangent sheaf $f^*TX$ is ample. Such rational curves are called \emph{very free} curves.
\end{definition}

We refer to Koll\'ar's book \cite{Kollar} or the work of Koll\'ar-Miyaoka-Mori \cite{KMM} for the background. If $X$ is separably rationally connected, then $X$ is rationally connected. The converse is true when the ground field is of characteristic zero by using the generic smoothness for the dominant map $u^{(2)}$. In positive characteristics, the converse statement is open.

In characteristic zero, a very important class of rationally connected varieties are Fano varieties, i.e., smooth varieties with ample anticanonical bundles. In positive characteristic, we only know that they are rationally chain connected. 

\begin{question}[Koll\'ar] In arbitrary characteristic, are Fano varieties separably rationally connected?
\end{question}

The question is open even for Fano hypersurfaces in projective spaces. In this paper, we prove the following theorem.

\begin{theorem}\label{mainn}
In arbitrary characteristic, a general Fano hypersurface of degree $n$ in $\P^n_k$ contains a minimal very free rational curve of degree $n$, i.e., the pullback of the tangent bundle has the splitting type $\O(2)\oplus\O(1)^{\oplus(n-2)}$.
\end{theorem}

\begin{theorem}\label{main}
In arbitrary characteristic, a general Fano hypersurface in $\mathbb{P}^n_k$ is separably rationally connected. 
\end{theorem}

de Jong-Starr \cite{dJS1} proved that every family of separably rationally connected varieties over a curve admits a rational section. Thus using Theorem \ref{main}, we give another proof of Tsen's theorem. 

\begin{cor}
Every family of Fano hypersurfaces in $\P^n$ over a curve admits a rational section.\qed
\end{cor}

\begin{ack} The author would like to thank his advisor Professor Jason Starr for helpful discussions. \end{ack}

\section{Typical Curves and Deformation Theory}

Let $n$ be an integer $\ge 3$. Let $X$ be a hypersurface of degree $n$ in $\mathbb{P}^n$. Let $C$ be a smoothly embedded rational curve of degree $e$ in $X$. We have the normal bundle short exact sequence.
$$\begin{CD}
0@>>>TC@>>>TX|_C@>>>\mathcal{N}_{C|X}@>>>0\\
\end{CD}$$
By adjunction, the degree of $TX|_C$ is the degree of $\O_{\P^n}(1)|_C$. Thus the degree of the normal bundle $\mathcal{N}_{C|X}$ is $e-2$ and the rank is $n-2$. 

\begin{definition}\label{typicaldef}
Let $e$ be a positive integer $\le n$. A smoothly embedded rational curve $C$ of degree $e$ in $X$ is \emph{typical}, if the normal bundle is the following:
\begin{equation*}\mathcal{N}_{C|X}\cong\big\{
\begin{array}{ll}
\mathscr{O}^{\oplus(n-3)}\oplus\mathscr{O}(-1),&\text{if } e= 1,\\
\mathscr{O}^{\oplus(n-e)}\oplus\mathscr{O}(1)^{\oplus(e-2)}, &\text{if } e\ge 2.
\end{array}\end{equation*}
The curve $C$ is a \emph{typical line} if the degree of $C$ is one. 
\end{definition}

Note that when $e=n$, typical rational curves of degree $n$ are very free. 

\begin{lemma}\label{typicalline}
Let $L$ be a smoothly embedded line in a hypersurface $X$ of degree $n$. Then $L$ is typical if and only if both of the following conditions hold:
\begin{enumerate}
\item $h^1(C, \mathcal{N}_{L|X})=0$,
\item $h^1(C,\mathcal{N}_{L|X}(-1))\le1$.
\end{enumerate}
\end{lemma}

\proof We may assume that $\N_{L|X}\cong \mathcal{O}(a_1)\oplus\dots\oplus\mathcal{O}(a_{n-2})$, where $a_1\ge\cdots\ge a_{n-2}$. Condition (1) is equivalent to that $a_1\ge\cdots\ge a_{n-2}\ge -1$. Together with condition (2), $a_{n-2}$ is either $0$ or $-1$. When $a_{n-2}=0$, $\N_{L|X}$ is semipositive, contradicting with the fact that the degree of $\N_{C|X}$ is $-1$. When $a_{n-2}=-1$, $\N_{L|X}/\O(a_{n-2})$ is semipositive. Because of the degree of the normal bundle, $L$ is typical.\qed

\begin{lemma}\label{typicalcurve}
Let $C$ be a smoothly embedded rational curve of degree $e$ in a hypersurface $X$ of degree $n$, where $2\le e\le n$. Then $C$ is typical if and only if both of the following conditions hold:
\begin{enumerate}
\item $h^1(C, \mathcal{N}_{C|X}(-1))=0$,
\item $h^1(C,\mathcal{N}_{C|X}(-2))\le n-e$.
\end{enumerate}
\end{lemma}

\proof 
Recall that the rank of the normal bundle $\mathcal{N}_{C|X}$ is $n-2$ and the degree is $e-2$. We may assume that $\mathcal{N}_{C|X}\cong \mathcal{O}(a_1)\oplus\dots\oplus\mathcal{O}(a_{n-2})$, where $a_1\ge\cdots\ge a_{n-2}$. Condition (1) is equivalent to that $a_{n-2}\ge 0$. Condition (2) implies that at most $n-e$ of $a_i$'s are $0$. By degree count, $C$ is a typical rational curve of degree $e$. \qed

Typical rational curves in the hypersurface $X$ are deformation open as very free curves in the following sense.

Let $H_n$ be the Hilbert scheme of hypersurfaces of degree $n$ in $\P^n$. It is isomorphic to some projective space. Let $\X\rightarrow H_n$ be the universal hypersurface. The morphism $\X\rightarrow H_n$ is flat projective and there exists a relative very ample invertible sheaf $\O_\X(1)$ on $\X$.

Let $R_{e,n}$ be the Hilbert scheme parameterizing flat projective families of one-dimensional subschemes in $\X$ with the Hilbert polynomial $P(d)=ed+1$. By \cite{Kollar} Theorem 1.4, $R_{e,n}$ is projective over $H_n$. 

Let $\C$ be the universal families over $R_{e,n}$, denoted by $\pi : \C\rightarrow R_{e,n}$. We have the following diagram,

$$
\xymatrix{
\C \ar[d]_\pi \ar[r]^<<<<<<i & R_{e,n}\times_{H_n} \X \ar[ld] \\
R_{e,n}&}
$$
where $i$ is a closed immersion.

\begin{prop}\label{typicaldefopen}
Let $e$ be a positive integer $\le n$. There exists an open subset in $R_{e,n}$ parameterizing typical curves of degree $e$ in hypersurfaces of degree $n$.
\end{prop}

\proof
Every typical curve of degree $e$ in a hypersurface of degree $n$ gives a point in $R_{e,n}$. Any small deformation of a smoothly embedded rational curve is still a smoothly embedded rational curve. Thus the proposition follows by Lemma \ref{typicalline}, Lemma \ref{typicalcurve} and the upper semicontinuity theorem \cite{Hartshorne} III.12.8.\qed 

\begin{lemma}\label{smoothdefopen}
There exists an open subset in $R_{e,n}$ such that for every closed point $(C,X)$ in the open subset, $C$ lies in the smooth locus of $X$.
\end{lemma}

\proof Let $S\subset \X$ be the relative singular locus in the universal hypersurface. $S$ is a closed subset of $\X$. Since $\pi$ is proper, the locus $\pi(i^{-1}(R_{e,n}\times_{H_n} S))$ is a closed subset of $R_{e,n}$ parametrizing the point $(C',X')$ such that $C'$ intersects the singular locus of $X'$. Thus the complement $U$ is open in $R_{e,n}$ and satisfies the desired property. \qed




Let $L$ be a typical line in a hypersurface $X$ of degree $n$ in $\P^n$. By definition, $\mathcal{N}_{L|X}\cong\mathscr{O}^{\oplus(n-3)}\oplus\mathscr{O}(-1)$. We have a canonically defined \emph{trivial subbundle} $\O^{\oplus(n-2)}$ of $\N_{L|X}$.


\begin{prop}\label{typicalconic}
Let $X$ be a hypersurface of degree $n$ in $\P^n$. Let $L$ and $M$ be two typical lines in $X$ intersecting transversally at only one point $p$. Assume that the following conditions hold:
\begin{enumerate}
\item the direction $T_p L$ is not in the trivial subbundle of $\mathcal{N}_{M|X}$;
\item the direction $T_p M$ is not in the trivial subbundle of $\N_{L|X}$.
\end{enumerate}
Then the pair $(L\cup M, X)\in R_{2,n}$ can be smoothed to a pair $(C, X^\prime)$ where $C$ is a typical conic in $X^\prime$. Furthermore, there exists an open neighborhood of $(L\cup M,X)$ in which any smoothing of $(L\cup M,X)$ is a typical conic.


\end{prop} 
\proof Let $D$ be the union of the lines $L$ and $M$. Since $D$ is a local complete intersection and lies in the smooth locus of $X$, the normal bundle $\mathcal{N}_{D|X}$ is locally free. We have the following short exact sequence.
$$\begin{CD}
0@>>>\mathcal{N}_{L|X}@>>>\mathcal{N}_{D|X}|_L@>>>T_p L\otimes T_p M@>>>0
\end{CD}$$
By \cite{GHS} Lemma 2.6, the locally free sheaf $\mathcal{N}_{D|X}|_L$ is the sheaf of rational sections of $\mathcal{N}_{L|X}$ which has at most one pole at the direction of $T_p M$. Since $\mathcal{N}_{L|X}\cong\mathscr{O}^{\oplus(n-3)}\oplus\mathscr{O}(-1)$, condition (2) implies that $\N_{D|X}|_L$ is isomorphic to $\O^{\oplus(n-2)}$. 

By the same argument, condition (1) implies that the sheaf $\mathcal{N}_{D|X}|_M$ is isomorphic to $\mathscr{O}^{\oplus(n-2)}$. Now we have the following short exact sequence. 
$$\begin{CD}
0@>>>\mathcal{N}_{D|X}|_M(-p)@>>>\mathcal{N}_{D|X}@>>>\mathcal{N}_{D|X}|_L@>>>0\\
@.@|@.@|@.\\
 @. \mathscr{O}(-1)^{\oplus(n-2)}  @. @.\mathscr{O}^{\oplus(n-2)}  @.
\end{CD}$$

First we claim that $D$ can be smoothed. Since $h^1(D,\N_{D|X})=0$, the pair $(D,X)$ is unobstructed in $R_{2,n}$, cf. \cite{Kollar} I.2. By \cite{Starr0} Lemma 3.17, it suffices to show that the map 
$$H^0(D,\mathcal{N}_{D|X})\rightarrow H^0(L,\N_{D|X}|L)\rightarrow T_p L\otimes T_p M$$
is surjective. Since $H^1(M,\mathcal{N}_{D|X}|_M(-p))=0$, the first map is surjective. Since $H^1(L,\mathcal{N}_{D|X}|_L)=0$, the second map is surjective.

Let $q$, $r$ be two distinct points on $L-\{p\}$. By the long exact sequence associated to the above short exact sequence at $h^1$, we get $h^1(D,\mathcal{N}_{D|X}(-q))=0$ and $h^1(D,\mathcal{N}_{D|X}(-q-r))=n-2$. 

Now for any smoothing $(D_t,X_t)$ of $(D,X)$ over $T$, we can specify two distinct points $p_t$ and $q_t$ on $D_t$ which specialize to $q$ and $r$ on $D$. By Lemma \ref{smoothdefopen}, after shrinking $T$, the conic $D_t$ lies in the smooth locus of $X_t$. Thus $D_t$ is smoothly embedded. By the upper semicontinuity theorem and Lemma \ref{typicalcurve}, $D_t$ is a typical conic in $X_t$. \qed

\begin{definition}\label{tc}
Let $X$ be a hypersurface of degree $n$ in $\P^n$. A \emph{typical comb} with $m$ teeth in $X$ is a reduced curve in $X$ with $m+1$ irreducible components $C, L_1,\cdots,L_m$ satisfying the following conditions:
\begin{enumerate}
\item $C$ is a typical conic in $X$;
\item $L_1,\cdots,L_m$ are disjoint typical lines in $X$ and each $L_i$ intersects $C$ transversally at $p_i$.
\end{enumerate}
The conic $C$ is called the \emph{handle} of the comb and $L_i$'s are called the \emph{teeth}.
\end{definition}

\begin{prop}\label{typicalcomb}
Let $X$ be a hypersurface of degree $n$ in $\P^n$. Let $D=C\cup L_1\cup\cdots\cup L_{n-2}$ be a typical comb with $n-2$ teeth in $X$. Let $p_i$ be the intersection point $L_i\cap C$. Assume that the following conditions hold:
\begin{enumerate}
\item the direction $T_{p_i} C$ is not in the trivial subbundle of $\mathcal{N}_{L_i|X}$;
\item the directions $T_{p_i} L_i$ are general in $\N_{C|X}$ such that the sheaf $\N_{D|X}|_C$ is isomorphic to $\O(1)^{\oplus(n-2)}$.
\end{enumerate}
Then the pair $(D, X)\in R_{n,n}$ can be smoothed to a pair $(C', X^\prime)$ where $C'$ is a very free curve in $X^\prime$. 
\end{prop} 

\proof The proof is very similar to the proof of Proposition \ref{typicalconic}. Here we only sketch the proof. Condition (1) implies that the sheaf $\N_{D|X}|_{L_i}$ is isomorphic to $\O^{\oplus(n-2)}$ for each $i$. We have the following short exact sequence.
$$\begin{CD}
0@>>>\sqcup_i\mathcal{N}_{D|X}|_{L_i}(-p)@>>>\mathcal{N}_{D|X}@>>>\mathcal{N}_{D|X}|_C@>>>0\\
@.@|@.@|@.\\
 @. \mathscr{O}(-1)^{\oplus(n-2)}  @. @.\mathscr{O}(1)^{\oplus(n-2)}  @.
\end{CD}$$
Since $H^1(D,\N_{D|X})=0$, $D$ is unobstructed. By diagram chasing, the map $H^0(D,\N_{D|X})\rightarrow \sqcup_i T_{p_i}C\otimes T_{p_i} L_i$ is surjective. Thus we can smooth the typical comb $D$.

Now we may choose a smoothing $(D_t,X_t)$ and specify two distinct points $(q_t, r_t)$ which specialize to two distinct points $(q,r)$ on $C-\{p_1,\cdots,p_{n-2}\}$. By the long exact sequence, we know that $h^1(D, \N_{D|X}(-q-r))=0$. By Lemma \ref{smoothdefopen} and the upper semicontinuity theorem, a general smoothing of the pair $(D,X)$ gives a very free curve in a general hypersurface.\qed

\section{An Example}
In this section, we construct a hypersurface of degree $n$ in $\P^n$, which contains a special configuration of lines. Later we will use this example to produce a very free curve in a general hypersurface.

Let $n$ be an integer $\ge 4$. Let $\lbrack x_0:\dots:x_n\rbrack$ be the homogeneous coordinates for $\mathbb{P}^n$. Let $X$ be a hypersurface of degree $n$ in the projective space $\P^n$ defined by the following equation.

\begin{eqnarray*}
\begin{array}{lllll}
x_0^{n-1}x_n 
&+x_1^{n-3}x_n^2x_0
& +(x_1^{n-1}+x_0x_1^{n-2}+\dots+x_0^{n-3}x_1^2)x_2 
&+(x_2^{n-1}+x_0x_2^{n-2}+\dots+x_0^{n-3}x_2^2)x_3
&+\dots\\

& +x_1^{n-4}x_n^3x_3
 &+(x_0x_1^{n-2}+\dots+x_0^{n-3}x_1^2)x_3
 &+(x_0x_2^{n-2}+\dots+x_0^{n-3}x_2^2)x_4&+\dots\\
&\vdots &\vdots\\

&+x_1x_n^{n-2}x_{n-2}
&+(x_0^{n-4}x_1^{3}+x_0^{n-3}x_1^2)x_{n-2} &+(x_0^{n-4}x_2^{3}+x_0^{n-3}x_2^2)x_{n-1}&+\dots\\

&+x_n^{n-1}x_{n-1}&+x_0^{n-3}x_1^2x_{n-1} &+x_0^{n-3}x_2^2x_{1}&+\dots
\end{array}\end{eqnarray*}

\begin{notation}\label{spiky}Let $p$ be the point $\lbrack 1:0:\dots:0\rbrack$ and $q$ be the point $\lbrack 0:1:0:\dots:0\rbrack$. Let $L_i$ be the line spanned by $\{e_0,e_i\}$ for $i=1,\dots, n-1$ and $L_n$ be the line spanned by $\{e_1,e_n\}$. It is easy to check that they all lie in the hypersurface $X$. Let $C$ be the union of $L_1,\cdots, L_n$. The following picture shows the configuration of the points and the lines in the projective space. \end{notation}

 \begin{center}
  \includegraphics[scale=0.5]{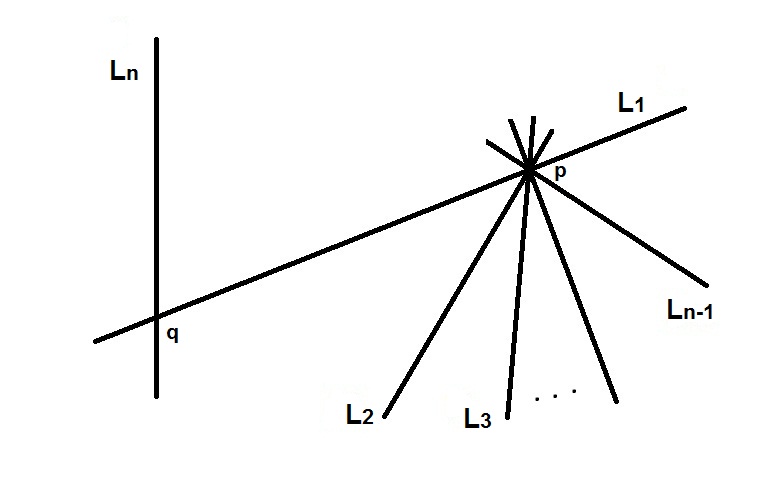}  

\end{center}

\begin{lemma} 
\begin{enumerate}
\item Both $p$ and $q$ lie in the smooth locus of $X$. 
\item The tangent space $T_p X$ is the hyperplane $\{x_n=0\}$, which is spanned by the lines $L_1,\dots,L_{n-1}$. 
\item The tangent space of $T_q X$ is the hyperplane $\{x_2=0\}$.
\end{enumerate}
\end{lemma}

\proof By taking the partial derivatives of $F$, we have 
$\frac{\partial F}{\partial x_i}(p)=0$ for $i=0,\cdots,n-1$ and $\frac{\partial F}{\partial x_n}(p)=1$. Similarly, we have $\frac{\partial F}{\partial x_i}(q)=0$ for $i\neq 2$ and $\frac{\partial F}{\partial x_2}(q)=1$.\qed


\begin{lemma}
The lines $L_1,\cdots, L_{n-1}$ are in the smooth locus of $X$. 
\end{lemma}

\proof We will prove the case for line $L_1$. The remaining cases can be computed directly by the same method. Denote $L_1=\{[x_0:x_1:0:\dots:0]\in\mathbb{P}^n\}$. By restricting the partial derivatives of the defining equation of the hypersurface $X$ on $L_1$, we get the following.
\begin{equation}
\begin{array}{l}
 \frac{\partial F}{\partial x_2}|_{L_1}=x_1^{n-1}+x_0x_1^{n-2}+\dots+x_0^{n-3}x_1^2
\\ \frac{\partial F}{\partial x_3}|_{L_1}=x_0x_1^{n-2}+\dots+x_0^{n-3}x_1^2
\\ \vdots\\ \frac{\partial F}{\partial x_{n-2}}|_{L_1}=x_0^{n-4}x_1^{3}+x_0^{n-3}x_1^2
\\  \frac{\partial F}{\partial x_{n-1}}|_{L_1}=x_0^{n-3}x_1^2
\\\frac{\partial F}{\partial x_n}|_{L_1}=x_0^{n-1}
\end{array}\label{L1}\end{equation}

For points on $L_1$ with $x_0\neq 0$, we have $\frac{\partial F}{\partial x_n}|_{L_1}\neq 0$. At the point $q$, $\frac{\partial F}{\partial x_2}|_{L_1}\neq 0$. Hence every point on the line $L_1$ is a smooth point of $X$. \qed

\begin{lemma}
The line $L_n$ is in the smooth locus of $X$. 
\end{lemma}

\proof By restricting the partial derivatives of the defining equation of $X$ on $L_n$, we get the following.
\begin{equation}
\begin{array}{l}
 \frac{\partial F}{\partial x_0}|_{L_n}=x_1^{n-3}x_n^2
\\ \frac{\partial F}{\partial x_3}|_{L_n}=x_1^{n-4}x_n^3
\\ \vdots\\ \frac{\partial F}{\partial x_{n-2}}|_{L_n}=x_1x_n^{n-2}
\\  \frac{\partial F}{\partial x_{n-1}}|_{L_n}=x_n^{n-1}
\\\frac{\partial F}{\partial x_2}|_{L_n}=x_1^{n-1}
\end{array}\label{Ln}\end{equation}

For points on $L_n$ with $x_1\neq 0$, we have $\frac{\partial F}{\partial x_2}|_{L_n}\neq 0$. For points on $L_n$ with $x_n\neq 0$, we have $\frac{\partial F}{\partial x_{n-1}}|_{L_n}\neq 0$.  Hence every point on the line $L_n$ is a smooth point of $X$. \qed



\begin{prop}\label{x0p1}With the notations as above, $X$ satisfies the following properties.
\begin{enumerate}
\item The lines $L_1,\dots,L_n$ are typical in $X$. 
\item For $i=1,\cdots, n-1$, the trivial subbundle of the normal bundle $\mathcal{N}_{L_i|X}$ at $p$ is generated by $\partial_{\overline{i+1}}-\partial_{\overline{i+2}},\dots,\partial_{\overline{i+n-3}}-\partial_{\overline{i+n-2}}$, where $\overline{j}$ takes value in $1,\dots, n-1$ mod $n-1$.
\item The trivial subbundle of the normal bundle $\mathcal{N}_{L_1|X}$ at $q$ is generated by $\partial_3,\dots, \partial_{n-1}$
\item The trivial subbundle of the normal bundle $\mathcal{N}_{L_n|X}$ at $q$ is generated by $\partial_3,\dots,\partial_{n-1}$.

\end{enumerate}
\end{prop}

\proof Let $L$ be a line in $X$. We have the following short exact sequences.
$$\begin{CD}
0 @>>> \mathcal{N}_{L|X}(-1) @>>> \mathcal{N}_{L|\mathbb{P}^n}(-1) @>>> \mathcal{N}_{X|\mathbb{P}^n}|_{L}(-1) @>>> 0\\
  @. @| @| @| @. \\
0 @>>>  \mathcal{N}_{L|X}(-1) @>>> \mathscr{O}_{L}^{\oplus(n-1)} @>>> \mathscr{O}_{L}(n-1)@>>> 0
\end{CD}$$
The associated long exact sequence is the following.  
$$\begin{CD}
 H^0(L,  \mathcal{N}_{L|X}(-1)) \rightarrow k^n@>\alpha>> H^0(L,  \mathscr{O}(n-1))\rightarrow H^1(L,  \mathcal{N}_{L|X}(-1))\rightarrow0
\end{CD}$$
where the map $\alpha$ sends the natural basis of $k^n$ to the derivatives of $F$ with respect to the normal directions of $L$ in $\P^n$. By Lemma \ref{typicalline}, $L$ is typical if and only if the image of $\alpha$ is of codimension one in $H^0(L,\mathscr{O}(n-1))$.

When $L=L_1$, by (\ref{L1}), $\frac{\partial F}{\partial x_2}|_{L_1},\dots, \frac{\partial F}{\partial x_n}|_{L_1}$ form a codimensional-one subspace of $H^0(L_1,\mathscr{O}_{L_1}(n-1))$. Thus we get that $H^1(L_1,\mathcal{N}_{L_1|X}(-1) )$ is one dimensional, i.e., $L_1$ is typical in $X$. 

By the short exact sequence above, $N_{{L_1|X}}(-1)$ is a subbundle of the trivial bundle  $\mathscr{O}_{L_1}^{\oplus(n-1)}$ which maps to $0$ in $\mathscr{O}_{L_1}(n-1)$. Let $\partial_2,\dots,\partial_n$ 
be the generators of $\mathscr{O}_{L_1}^{\oplus(n-1)}$. We get $N_{{L_1|X}}(-1)$ is generated by $x_0(\partial_2-\partial_3)-x_1(\partial_3-\partial_4), \dots, x_0(\partial_{n-2}-\partial_{n-1})-x_1\partial_{n-1}, x_0^2\partial_{n-1}-x_1^2\partial_n$ as an $\O_{L_1}$-module. If we restrict the bundle at $p$ and $q$, we get (2) and (3) for $L_1$. 

When $L=L_2,\cdots, L_{n-1}$, we can prove in a similar way. When $L=L_n$, (4) follows from the same computation as above by applying (\ref{Ln}).\qed  


With the description of the trivial subbundles of the normal bundles of lines in $X$ as above, we get the following corollaries.

\begin{cor}\label{cond1} We have the following statements.
\begin{enumerate}
\item The lines $L_1$ and $L_n$ are typical in $X$.
\item The direction $T_q L_1$ is not in the trivial subbundle of $\N_{L_n|X}$.
\item The direction $T_q L_n$ is not in the trivial subbundle of $\N_{L_1|X}$.\qed
\end{enumerate}
\end{cor}

\begin{cor} \label{cond2}We have the following statements.
\begin{enumerate}
\item The lines $L_2,\cdots,L_{n-2}$ are typical in $X$.
\item The direction $T_p L_1$ is not in the trivial subbundle of $\N_{L_i|X}$ for $2\le i\le n-1$.
\item The directions $T_p L_2,\cdots,T_p L_{n-1}$ span the normal bundle $\N_{L_1|X}$ at $p$.\qed
\end{enumerate}
\end{cor}
%


\section{Proof of the Main Theorem}
\begin{lemma}[\cite{Harris} Ex 13.8]\label{p=0} Let $C$ be the union of $n$ lines $L_1,\cdots,L_n$ in $\P^n$ as in Notation \ref{spiky}. Then the Hilbert polynomial of $C$ is $P(d)=nd+1$. In particular, the arithmetic genus of $C$ is $0$.
\end{lemma}

\proof This can be computed directly. For any $d>0$, when $i=1,\cdots,n-1$, the homogeneous polynomials of degree $d$ that do not vanish on $L_i$ are generated by $\{x_0^d, x_0^{d-1}x_i,\cdots,x_i^d\}$. The homogeneous polynomials of degree $d$ that do not vanish on $L_n$ are generated by $\{x_1^d, x_1^{d-1}x_i,\cdots,x_n^d\}$. Thus when $d\gg 0$, $P(d)=H^0(C,\O_C(d))=nd+1$.\qed

The curve $C$ is an example of curves with rational $n$-fold point, cf. \cite{dawei} 3.7. The following lemma is an analogue of \cite{dawei} Lemma 3.8.

\begin{lemma}\label{h0C} With the same notations as in Lemma \ref{p=0}, the following properties hold for $C$ for every positive integer $d$:
\begin{enumerate}
\item $h^0(C, \O_C(d))=nd+1$ and $h^1(C,\O_C(d))=0$.
\item $h^1(C,\I_C(d))=0$.
\item $h^0(C,\I_C(d))=h^0(\P^n,\O(d))-nd-1$. 
\end{enumerate}
\end{lemma}

\proof By Riemann-Roch and Lemma \ref{p=0}, we have $$h^0(C,\O_C(d))\ge \chi(\O_C(d))=nd+1-p_g(C)=nd+1.$$ 
On the other hand, every global section of $\O_C(d)$ is obtained by gluing global sections on each component, which imposes at least $n-1$ linear conditions. Since we have $h^0(L_i,\O(d))=d+1$ for every $i$, $$h^0(\O_C(d))\le n(d+1)-(n-1)=nd+1.$$ 
Therefore, $h^0(C, \O_C(d))=nd+1$ for every positive integer $d$. Since the image of $r:H^0(\P^n,\O(d))\rightarrow H^0(C,\O_C(d))$ has dimension $nd+1$ as in Lemma \ref{p=0}, the map $r$ is surjective. The lemma follows by considering the long exact sequence associated to the ideal sheaf of $C$. \qed

\begin{construction}
Let $C$ be the union of $n$ lines $L_1,\cdots,L_n$ in $\P^n$ as in Notation \ref{spiky}. If we consider $L_1\cup L_n$ as a conic in $\P^n$, there exists a smooth affine pointed curve $(T,0)$ and a smoothing  $D'\rightarrow (T,0)$ satisfying the following conditions:
\begin{enumerate}
\item The special fiber $D'_0$ is $L_1\cup L_n$;
\item For any $t\in T-\{0\}$, $D'_t$ is a smooth conic contained in the plane spanned by $L_1$ and $L_n$.
\end{enumerate}
We can choose $n-2$ sections $s_i:(T,0)\rightarrow D'$ for $i=1,\cdots, n-2$ such that $s_i(0)=p$ for all $i$'s and for $t\in T-\{0\} $, $s_i(t)$'s are all distinct on $D'_t$. 

For any $s_i(t)$, there exists a unique line $L_{i+1}(t)$ through $s_i(t)$ parallel to $L_{i+1}$. After gluing the families of lines $L_{i+1}(t)$ on $D'_t$ at $s_i(t)$ for all $i$'s, we get a family of reducible curves $\pi:D\rightarrow (T,0)$ satisfying the following conditions:
\begin{enumerate}
\item The special fiber $D_0$ is $C$ as constructed in \ref{spiky}.
\item For any $t\in T-\{0\}$, $D_t$ is a comb with the handle $D'_t$ and with the teeth lines.
\end{enumerate} 

The family $\pi:D\rightarrow (T,0)$ is flat by Lemma \ref{p=0}. We have the following diagram.

$$\xymatrix{D_0=C \ar[r] \ar[d] & D \ar[d]_\pi \ar[r]^i  &\P^n_T  \ar@{->}[ld]^\pi  \\0 \ar[r] &(T,0)}$$
\end{construction}

\begin{lemma}\label{flatlift} Let $\I_D$ be the ideal sheaves of $D$ in $\P^n_T$. The sheaf $\pi_*\I_D(d)$ is locally free on $T$ for any $d>0$. 
\end{lemma}

\proof By the cohomology and base change theorem \cite{Hartshorne} III.12.9, it suffices to show that $h^0(\P^n_t, I_{D_t}(d))$ is constant. By upper semicontinuity and Lemma \ref{h0C}, we have $h^0(\P^n_t, I_{D_t}(d))\le h^0(\P^n,\O(d))-nd-1$. On the other hand, for any $t\in T-\{0\}$, the curve $D_t$ is a local complete intersection and $\O_{D_t}(d)$ is a positive bundle on $D_t$. Thus we have $h^1(D_t,\O_{D_t}(d))=0$ and $h^0(D_t,\O_{D_t}(d))=nd+1$. Consider the following exact sequence.
$$\begin{CD}
0@>>>H^0(\P^n_t,\I_{D_t}(d))@>>>H^0(\P^n_t, \O(d))@>>>H^0(D_t,\O_{D_t}(d))
\end{CD}$$
We get $h^0(\P^n_t,\I_{D_t}(d))\ge h^0(\P^n,\O(d))-nd-1$.\qed

\proof [Proof of Theorem \ref{mainn}]
The theorem is trivial for $n=2,3$. We can assume that $n\ge 4$. By \cite{Kollar} IV.3.11 and Lemma \ref{smoothdefopen}, it suffices to produce one very free curve in a hypersurface of degree $n$. By Lemma \ref{flatlift}, after shrinking $T$, hypersurfaces of degree $n$ containing $D_t$ in $\P^n_t$ form a trivial projective bundle over $(T,0)$. Thus the family $\pi:D\rightarrow (T,0)$ admits a lifting to a flat family of pairs $\pi:(D,\X_T)\rightarrow (T,0)$ in $R_{n,n}$ such that the special fiber $(D_0,\X_0)$ is $(C,X)$ which is constructed in Section 3.
$$\xymatrix{ D \ar[d]_\pi \ar[r]^i  &\X_T \ar[r] \ar[ld] &\P^n_T \ar[lld]^\pi \\(T,0)}$$
All the following steps of the proof requires to shrink $T$ if necessary. By Proposition \ref{typicalconic} and Corollary \ref{cond1}, we may assume that the handle $D'_t$ is a typical conic in $\X_t$ for every $t\in T-\{0\}$. By Proposition \ref{typicaldefopen} and Corollary \ref{cond2} (1), all the teeth of the comb $D_t$ are typical. Thus for every $t\in T-\{0\}$, we get a typical comb $D_t$ as in Definition \ref{tc}. Now the theorem follows if we verify the two conditions in Proposition \ref{typicalcurve}. Since they are open conditions, it suffices to check on the special fiber $(C,X)$, which is proved in Corollary \ref{cond2}.\qed

\proof [Proof of Theorem \ref{main}]
 By \cite{Kollar} IV.3.11 and Lemma \ref{smoothdefopen}, it suffices to produce one very free curve in a hypersurface of degree $d$. Let $Y$ be a general smooth Fano hypersurface of degree $d$ in $\P^n$. When $d=n$, this is proved in Theorem \ref{mainn}. When $d<n$, 
we may choose a general linear subspace $L$ of dimension $d$ such that $Y\cap L$ is smooth and contains a very free curve $f:\P^1\rightarrow Y\cap L$ by Theorem \ref{mainn}. By the normal bundle exact sequence, 
$$\begin{CD}
0@>>>T(Y\cap L)@>>>TY@>>>\N_{Y\cap L|Y}@>>>0
\end{CD}$$
the sheaf $f^*T(Y\cap L)$ is positive and the sheaf $\N_{Y\cap L|Y}$ is isomorphic to $\N_{L|\P^n}$, which is $\O(1)^{\oplus(n-d)}$. Therefore the pullback bundle $f^*TY$ is positive. Thus $f:\P^1\rightarrow Y\cap L\rightarrow Y$ is a very free curve in $Y$.\qed

\bibliography{myref}	
\bibliographystyle{alpha}	

\end{document}